\documentclass{amsart}[11pt]

\usepackage[all]{xy}
\usepackage{amsmath}
\usepackage{amssymb}
\usepackage{MnSymbol}

\newtheorem{theorem}{Theorem}[section]
\newtheorem{lemma}[theorem]{Lemma}
\newtheorem{corollary}[theorem]{Corollary}
\theoremstyle{definition}
\newtheorem{definition}[theorem]{Definition}

\newtheorem{proposition}[theorem]{Proposition}

\newtheorem{manualtheoreminner}{Theorem}
\newenvironment{manualtheorem}[1]{%
  \renewcommand\themanualtheoreminner{#1}%
  \manualtheoreminner
}{\endmanualtheoreminner}

\theoremstyle{remark}

\usepackage{graphicx}
\graphicspath{ {./images/} }

\numberwithin{equation}{section}

\begin{document}
	
\title{Dehn functions for mapping tori of amalgam of free groups}

\author{Qianwen Sun}
\address{Department of Mathematics, Statistics, and Computer Science,
	University of Illinois at Chicago,
	322 Science and Engineering Offices (M/C 249),
	851 S. Morgan St.
	Chicago, IL 60607-7045
	 }
\email{qsun24@uic.edu}

\date{}
\maketitle

\textbf{Abstract}: For an amalgam of two free groups and a particular kind of automorphism, we show that the Dehn function of the corresponding mapping torus is quadratic.

\section{Introduction}

Dehn functions have attracted a lot of attention in recent years. Named after Max Dehn, a Dehn function is an optimal function that bounds the area of an identity in terms of the defining relators in a finitely presented group. It is also closely connected with algorithmic complexity of the word problem, which is the problem of deciding whether a given word equals to 1. For a finitely presented group, it has solvable word problem if and only if it has recursive Dehn function \cite{SM}.

A milestone that related Dehn functions and hyperbolic groups, proved by Gromov \cite{M}, is that a finitely presented group is hyperbolic if and only if its Dehn function is linear, or equivalently, it satisfies a linear isoperimetric inequality.

In the same paper, Gromov also pointed out that there exists an isoperimetric gap: if a finitely presented group satisfies a subquadratic isoperimetric inequality, then its Dehn function is linear. In particular, no group has a Dehn function equivalent to $n^d$ with $d\in(1,2)$.

Brinkmann \cite{B} classified all hyperbolic mapping tori of free groups. In 2000, Macura \cite{N} proved that the mapping torus of a polynomially growing automorphism of finitely generated free group satisfies a quadratic isoperimetric inequality.
Bridson and Groves \cite{MD} generalized it to an arbitrary automorphism in 2010. 

In this paper, following Bridson and Groves, we consider an amalgam of two free groups but only with a certain kind of automorphism.

Let $G$ be an amalgam of two finitely generated free groups $F_1$ and $F_2$ amalgamated along finitely generated subgroups $A_i\in F_i$, where $A_1$ is a malnormal subgroup. Let $\psi$ be an automorphism of $F_2$ that fixes $A_2$ pointwise. Then there exists an induced automorphism $\varphi$ of $G$ such that it restricts to identity on $F_1$, and $\psi$ on $F_2$.

\begin{theorem}\label{1}
The Dehn function of the mapping torus $$M_{\varphi}=\{G,t;t^{-1}at=\varphi (a), a\in G\}$$ is quadratic.
\end{theorem}

Note that when $A_2$ is generated by one element $u$ and $\psi=Ad_u$, the induced map $\varphi$ is the Dehn twist along the edge.

The group $G$ in Theorem~\ref{1} is hyperbolic. This is a special case of a more general result, see \cite{OA}, where they show that amalgam of two hyperbolic groups are hyperbolic if at least one of the amalgamated subgroup is malnormal, by proving such an amalgam satisfies a linear isoperimetric inequality. For our work, we show a stronger result: Single Bounded Theorem.  

\begin{theorem}\label{SBT}
Let $G$ be an amalgam of two free groups $F_1$ and $F_2$ amalgamated along finitely generated subgroups $A_1\in F_1$ and $A_2\in F_2$. If $A_1$ is a malnormal subgroup, then there exists a presentation of $G$ and a constant $B$ such that for any word $z=X_1Y_1\dots X_nY_n=_G1$, $X_i\in F_1$, $Y_i\in F_2$, ${\rm{Area}}(z)\leq B\sum^n_{i=1}|X_i|$ in this presentation.
\end{theorem}

Here, $|X_i|$ denotes the length of the reduced word in $F_1$ representing $X_i$. We make this convention throughout the rest of this paper.

As a corollary, we recover the above mentioned fact that $G$ is hyperbolic.

\begin{corollary}
Let $G$ be an amalgam of two free groups $F_1$ and $F_2$ amalgamated along finitely generated subgroups $A_1\in F_1$ and $A_2\in F_2$. If one of $A_i$ is a malnormal subgroup, then $G$ is hyperbolic.
\end{corollary}

Mapping tori have been powerful in low-dimensional topology. Since any outer automorphism $[\bar{f}]$ of a surface group is realized by a homeomorphism $f$ of the underlying surface. The mapping tori is fundamental group of a 3-manifold.

A deep result of Thurston \cite{T} states that in this case the 3-manifold $M_f$ is hyperbolic if and only if $f$ is a pseudo-Anosov homeomorphism of $S$. In all other cases, $\pi_1(M_f)$ has quadratic Dehn functions, for example, because it is automatic \cite{D}.

This paper is organized as follows: first in Section \ref{preliminary}, we recall some preliminaries. Then in Section~\ref{combinatorial}, we introduce some combinatorial methods and use them to prove the Single Bounded Theorem. Finally in Section \ref{proof}, we introduce some geometric methods and use them to prove Theorem~\ref{1}.

\section{Preliminaries}\label{preliminary}
In this section, we will recall some preliminary knowledge. Some standard references on these topics are \cite{MA}, \cite{R}.
\subsection{N-reduced set}

An N-reduced set is an important concept in the study of subgroups of a free group. We recall the definition and some properties here, see \cite{R}.

Consider a set of reduced words $W=\{w_1,w_2,\dots\}$ in a free group $F$ with a fixed basis. We call $W$ \emph{N-reduced} if for all triples $v_1,v_2,v_3$ of the form $w_i^{\pm1}$, the following conditions hold:

(N0) $v_1\neq 1$;

(N1) $v_1v_2\neq 1$ implies $|v_1v_2|\geq \max\{|v_1|, |v_2|\}$;

(N2) $v_1v_2\neq 1$ and $v_2v_3\neq 1$ implies $|v_1v_2v_3|>|v_1|-|v_2|+|v_3|$.

By establishing the following proposition, Nielsen \cite{N1} shows that every finite generated subgroup of a free group is free. 
\begin{proposition}\label{Nielsen}
	If $W=\{w_1,\dots, w_m\}$ is a finite set in a free group $F$, then $W$ can be carried by Nielsen transformations into some $W'$ such that $W'$ is N-reduced and $\langle W\rangle=\langle W'\rangle$
\end{proposition}

There is a slightly stronger result than Nielsen's original work due to Zieschang \cite{Z}.

\begin{proposition}\label{uncanceled}
If $W=\{w_1,w_2,\dots,w_m\}$ is N-reduced, then one may associate with each $v$ in $W^{\pm1}$ words $a(v)$ and $m(v)$, with $m(v)\neq 1$, such that
$v=a(v)m(v)a(v^{-1})^{-1}$ reduced,
and such that if
$u=v_1v_2\dots v_t$, $v_i \in W^{\pm1}$ and all $v_iv_{i+1}\neq 1$,
then $m(v_1),\dots, m(v_t)$ remain uncanceled in the reduced form of $u$.
\end{proposition}

\subsection{Amalgams}

Let $G_1$ and $G_2$ be groups. Let $i_1: H\hookrightarrow G_1$ and $i_2: H\hookrightarrow G_2$ be embeddings. Then the \emph{amalgamated free product of $G_1$ and $G_2$ over $H$} is the group 
$$G_1*_HG_2:=\langle G_1,G_2;i_1(a)=i_2(a),a\in H \rangle.$$

There are extensive studies on the structure of amalgams. One fact is that the two natural inclusions $G_1\rightarrow G_1*_HG_2$ and $G_2\rightarrow G_1*_HG_2$ are both embeddings. So $G_1$ and $G_2$ can be identified as subgroups of $G_1*_HG_2$.

Another classical result is the following lemma.

\begin{lemma}\label{2}[\cite{R}, IV. Theorem 2.6]
Let $G=G_1*_HG_2$. Let $X_1,\dots\dots,X_n\in G_1$ and $Y_1,\dots\dots,Y_n\in G_2$ be elements such that $X_i\notin H$ if $i>1$ and $Y_i\notin H$ if $i<n$ where $n>1$. Then $X_1Y_1\dots\dots X_nY_n\neq 1$ in $G$.
\end{lemma}
As a corollary, 
\begin{corollary}\label{R1}
Let $G=G_1*_HG_2$. Let $X_1,\dots\dots,X_n\in G_1$ and $Y_1,\dots\dots,Y_n\in G_2$ such that $z=X_1Y_1\dots\dots X_nY_n$ represents $1$ in $G$, then there exists $k$ such that $X_k\in H$ or $Y_k\in H$. Furthermore, if $n> 1$ and all $X_i$'s and $Y_i$'s are nontrivial except $X_1$ and $Y_n$, then there exists nontrivial $X_k$ or $Y_k$ that is in $H$.
\end{corollary}

\subsection{Dehn functions}

For a finitely presented group $G=\langle X;R\rangle$, a word $w\in F(X)$ represents $1$ in G if and only if it can be represented as finite product of conjugates of elements in $R$. Among all such representations, there is a smallest integer $n$ such that $w$ can be represented by a product of conjugates of $n$ elements in $R$. We call $n$ the \emph{area} of $w$. 

Van Kampen diagrams are a standard tool for the study of isoperimetric functions in groups, see \cite{M1} for details. The idea of a van Kampen Diagram is that for a word $w=_G1$, we have a simply connected planer 2-complex, whose 2-cells are labeled by relators and whose boundary is labeled by $w$. 

We call the number of 2-cells in a van Kampen diagram its \emph{area}. The area of $w$ defined above is the minimum area of a van Kampen diagram with boundary labeled by $w$.
 
For a finitely presented group $G=\langle X;R\rangle$, the \emph{Dehn function} is defined as $$D(n)=\max\{{\rm{Area}}(w)|w=_G1,|w|\leq n\}.$$ 
Although it depends on the particular presentation, different presentations of the same group give equivalent Dehn functions, see \cite{M1}. Here, two functions $f$ and $g$ are equivalent if $f\preceq g$ and $g\preceq f$, where $\preceq$ means that there exists a constant $C>0$ such that $$f(n)\leq Cg(Cn+C)+Cn+C.$$
The Dehn function of a group is defined as the equivalence class of functions.

\subsection{Mapping torus}
 Let us recall the definition of the mapping torus. For any automorphism $\varphi$ of a group $G$, the algebraic mapping torus is defined by 
 $$M_{\varphi}=\langle G,t;t^{-1}at=\varphi (a), a\in G\rangle.$$ Note that $G$ and $\langle t\rangle$ are embedded subgroups in the mapping torus.

In 2010, Bridson and Groves \cite{MD} proved the following theorem about mapping tori of free groups. 

\begin{theorem}\label{3}
If $F$ is a finitely generated free group and $\varphi$ is an automorphism of $F$, then the mapping torus $M_{\varphi}$ satisfies a quadratic isoperimetric inequality.
\end{theorem}

\section{Combinatorial methods}\label{combinatorial}

\subsection{Relating Pairs}

An N-reduced set $W=\{w_1,\dots, w_m\}$ in a free group $F$ forms a basis for the subgroup $\langle w_1,\dots, w_m\rangle$. We assume all $w_i$ are reduced words in this section. 

By a \emph{combination} we mean a word in $\langle W\rangle$ that is represented by a product of elements in $W^{\pm1}$. If no two adjacent elements are inverse of each other, the combination is said to be \emph{reduced}, in which case it is the unique representation in terms of this basis.

For a reduced combination $v=v_1\dots v_t$ where $v_i\in W^{\pm1}$, there may be cancellation between $v_i$ and $v_{i+1}$. By Proposition \ref{uncanceled}, there is a unique partition for each $v_i$: $v_i=a_{i-1}^{-1}s_ia_i$, where $a_0$ and $a_t$ are trivial, $s_i$ is nontrivial and contains $m(v_i)$, and $v_1\dots v_t=s_1\dots s_t$ with no cancellation between $s_i$ and $s_{i+1}$. We call each $s_i$
 a \emph{stem}, and each $a_i$ a \emph{twig} in this combination, where $a_{i-1}$ and $a_i$ are left and right twig of $v_i$ respectively. 
 
 If a word is a stem or a twig in some combination, we call it a stem or a twig in $W$.
 
A reduced combination, when expressed as a product of elements in $W^{\pm1}$, can be represented by a path in the Cayley graph of $F$. After Stalling foldings, the folded path with labels is called \emph{the graph of the combination}. We make the convention that all twigs are placed vertically, and all stems horizontally in the graph.
 
\begin{figure}
	\includegraphics[width=90mm]{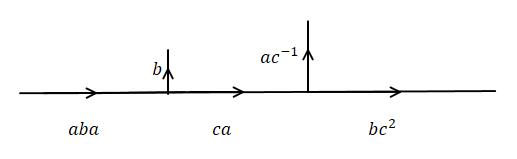}
	\caption{}
	\label{product}
\end{figure}
 
Figure \ref{product} is an example of such a graph in a free group generated by $a, b, c$, where $v_1=abab$, $v_2=b^{-1}ca^2c^{-1}$, $v_3=ca^{-1}bc^2$. Note that the horizontal line represents the reduced form of $v_1v_2v_3$.
 
For $v\in W^{\pm1}$, if it has a partition $asb$ in some combination, we call it a partition of $v$ in $W^{\pm1}$, and $s$ a stem of $v$.

\begin{definition}
A \emph{relating pair} is a pair of elements in $W^{\pm1}$: $[v_1,v_2]$, equipped with the following information:
\begin{enumerate}
\item
a partition in $W^{\pm1}$ for each $v_i$: $v_i=a_is_ib_i$;
\item
a nontrivial word $s$, which is called the \emph{relating segment}, that is an initial subword of $s_1$, and is also a subword of $s_2$. That is, $s_1=ss'$ and $s_2=lsm$ as strings of letters. Note that $s',l,m$ may be trivial.
\end{enumerate}

 A relating pair is \emph{trivial} if $v_1=v_2$, their partitions are the same, and $s=s_1=s_2$.
\end{definition}
 
The following lemma is obvious and we omit the proof here.

\begin{lemma}
If $W$ is finite, then there are only finitely many stems and twigs, and also finitely many relating pairs.
\end{lemma}

For the rest of this section, let $S$ denote the longest length of all stems, $T$ the longest length of all twigs, and $R$ the number of nontrivial relating pairs. 

\subsection{Identical Segments}

Let $z=z_1\dots z_n$ and $v=v_1\dots v_t$ be two reduced combinations. There exist stems $r_i$ of $z_i$, $s_i$ of $v_i$ in corresponding combinations, such that $z=r_1\dots r_n$ and $v=s_1\dots s_t$ are reduced form of $z$ and $v$. \emph{A segment of $z$} is a subword of $r_1\dots r_n$. If a segment of $z$ and a segment of $v$ are the same as words, we call them \emph{identical segments} between these two combinations. 

One can show identical segments by a graph. It consists of two graphs of combinations, one for $z$ on the top and one for $v$ on the bottom, while the corresponding labels in identical segments are in the same horizontal positions.

Identical segments from two ordered combinations generate a series of relating pairs as follows. First in $z$, suppose the stems that lie entirely in this segment are $r_k,\dots,r_{k+l}$. Consider the first letter of each stem, it corresponds to one letter in one stem of $v$, suppose these stems are $s_{i_k},\dots,s_{i_{k+l}}$ respectively. Note that some of them may not be entirely in this identical segment, and some adjacent stems may be the same one. For each $j$, $k\leq j\leq k+l$, there is a relating pair $[z_j,v_{i_j}]$, whose partition is the corresponding partition in each combinations, and the relating segment is the maximal initial subword of $r_j$ whose corresponding segment is in $s_{i_j}$, while this corresponding segment is the subword of $s_{i_j}$ which is the relating segment.

Note that this series may be empty if the segment in the first combination does not contain a single full stem. 

Call two identical segments \textit{matchless} if there is no trivial relating pair in this series. This definition does not depend on the order of two combinations.

\begin{lemma}\label{C}
	Let $W=\{w_1,\dots, w_m\}$ be an $N$-reduced set that generates a malnormal subgroup of the free group $F$. There exists a constant $C$ such that two combinations cannot have matchless identical segments that are longer than $C$.	
\end{lemma}

\begin{proof}
	
	Let $C=(R+3)S$. Suppose there are two combinations, say $x_1$ and $x_2$, that have mathchless identical segments $y$ of length longer than $C$. 
	
	Consider all the stems in $x_1$ that fully lie in $y$, the number of these stems is at least $(C-2S)/S\geq R+1$.
	That means there is a series of nontrivial relating pairs $p_1, p_2, \dots, p_n$ with $n\geq R+1$. So at least two of them are the same. We may assume that $p_{k+1}=p_1$ and this part is shown in Figure \ref{lines}, where $z_i$'s, $v_i$'s are the corresponding elements in $W^{\pm1}$.
	
	\begin{figure}
		\includegraphics[width=120mm]{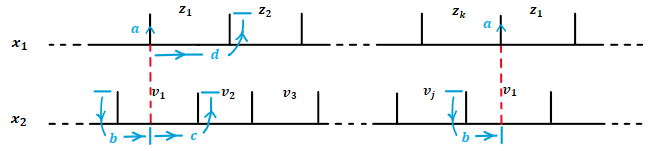}
		\caption{}
		\label{lines}
	\end{figure}
	
	Suppose $z_1=a^{-1}d$, $v_1=bc$ where $d$ and $c$ start with the relating segment of $p_1$. The part between two red lines is represented by both top and bottom segments so
	\begin{equation}\label{malnormal}
	\begin{split}
	&az_1\dots z_ka^{-1}=b^{-1}v_1\dots v_jb\\
	&\implies z_1\dots z_k=(ba)^{-1}v_1\dots v_jba.
	\end{split}
	\end{equation}
	
	First, $ba$ is not trivial. If $ba=1$, from \eqref{malnormal}, we get $z_1\dots z_k=v_1\dots v_j$. Since $W$ is a basis, we have $k=j$ and $z_i=v_i$ for all $i$.  
	
 Suppose $b$ ends with $l$, where $l$ is the initial part of the stem of $v_1$. Then the stem of $z_k$ also ends with $l$. Since $a=b^{-1}$, the right twig of $z_k$ starts with $l^{-1}$. If $l$ is not empty, the label of $z_k$ is not reduced, which is impossible because elements in $W$ are all in reduced form. Hence $l$ is empty and $b$ is the left twig of $v_1$. By $z_2=v_2$, the right twig of $z_1$ and $v_1$ are also the same. It follows that the relating pair $[z_1,v_1]$ is trivial, contradicting our assumptions.

	Next, $ba\notin \langle w_1,\dots, w_m\rangle$. We prove it by contradiction. Let $\overline{ba}$ represent the reduced form of $ba$. Suppose $\overline{ba}=y_1\dots y_l$ is a reduced combination, then $l\geq 1$ and
	
	\begin{equation}\label{cancelout0}
	\begin{split}
	&v_1^{-1}(y_1\dots y_l)z_1\\
	&=c^{-1}b^{-1}(\overline{ba})a^{-1}d\\
	&=\overline{c^{-1}d}.
	\end{split}
	\end{equation}

    From this reduction process, $y_1\dots y_l$ is freely reduced to $\overline{ba}$ and then this part is canceled completely. If $y_1\neq v_1$ and $y_l\neq z_1^{-1}$, we are done since that all $y_i$'s are canceled contradicts Proposition \ref{uncanceled}. Otherwise, there are 3 cases:
	
	\textbf{Case 1}: $y_1=v_1$ and $y_l\neq z_1^{-1}$. Then $ba=v_1\dots y_l\implies y_2\dots y_l=c^{-1}a$. Here $c^{-1}a$ must be a reduced word, otherwise the label of $z_1$ is not reduced. It also implies that $c^{-1}a$ is nontrivial because $c$ contains the relating segment, so $l\geq 2$.

	Now consider the combination $v_1y_2\dots y_lz_1$, in which no two adjacent elements are inverse of each other. To reduce it, we get
	\begin{equation}\label{cancelout}
	\begin{split}
	&v_1(y_2\dots y_l)z_1\\
	&=bc(c^{-1}a)a^{-1}d\\
	&=bd,
	\end{split}
	\end{equation}
where $y_2\dots y_l$ are all canceled. This contradicts Proposition \ref{uncanceled}.
	
	\textbf{Case 2}: $y_1\neq v_1$ and $y_l=z_1^{-1}$. Follow a similar proof as that in case 1, one can show that $y_1\dots y_{l-1}$ are all canceled in $v_1^{-1}y_1\dots y_{l-1}z_1^{-1}$. We omit the proof here.
	
    \textbf{Case 3}: $y_1=v_1$ and $y_l=z_1^{-1}$. Then $ba=v_1y_2\dots y_{l-1}z_1^{-1}\implies y_2\dots y_{l-1}=\overline{c^{-1}d}$.
    
If $c^{-1}d=1$, $z_2\dots z_kz_1=v_2\dots v_jv_1$. Then follow a similar analysis as that for the $ba$, it will result in a trivial pair. 

If $c^{-1}d$ is not trivial, then $l\geq 3$. Now in the following reduction process:
	\begin{equation}\label{cancelout3}
	\begin{split}
	&v_1(y_2\dots y_{l-1})z_1^{-1}\\
	&=bc(\overline{c^{-1}d})d^{-1}a\\
	&=\overline{ba},
	\end{split}
	\end{equation}
 $y_2\dots y_{l-1}$ are all canceled, which contradicts Proposition \ref{uncanceled}.

	Now $ba\notin \langle W\rangle$. Equation \eqref{malnormal} implies that $\langle W\rangle$ is not a malnormal subgroup, which contradicts our assumption. So two combinations cannot have matchless identical segments that are longer than $C$.	
\end{proof}

\subsection{Factors in An Amalgam}

Following the notation from last section, we make the following definition.

\begin{definition}
A word $X\in F$ is \emph{$W$-free} if $X$ has the shortest length among all elements in $\langle W \rangle X \langle W \rangle$.
\end{definition}

\begin{lemma}\label{D}
	Let $G=\langle F_1*G_1;w_i=u_i,i=1,\dots, m \rangle$ be an amalgam of a free group $F_1$ and a group $G_1$, $w_i\in F_1$, $u_i\in G_1$, and $W=\{w_i\}$ is an N-reduced basis that generates a malnormal subgroup of $F_1$. Then there exists a constant $D$ such that for any nontrivial word $X_1Y_1\dots X_nY_n=_G1$ with all $X_i$ $W$-free, and all $X_i$, $Y_i$ nontrivial except $X_1$ and $Y_n$, there exists $k$ such that $Y_k=_G\Pi u_{t_i}^{s_i}$ with $0<\sum|s_i|\leq D$.
\end{lemma}

\begin{proof}

We prove that it suffices to take $D=3C+2S+2T$, where $C$ is the constant in Lemma \ref{C}.

If $n=1$, by the fact $X_1$ is $w$-free, nontrivial and Corollary \ref{R1}, $Y_1\in\langle U\rangle$, where $U=\{u_1, \dots, u_m\}$. Change each $u_i$ to $w_i$, it follows that $X_1 W_1=1$ for some $W_1\in\langle W\rangle$. So $X_1=Y_1=1$, contradicting the assumption that $X_1Y_1$ it a nontrivial word.

So $n\geq2$. By Corollary \ref{R1}, there exists $k$ such that $X_k\in \langle W\rangle$ is nontrivial or $Y_k\in \langle U\rangle$ is nontrivial. If some $X_i$ is both $W$-free and in $\langle W\rangle$, then it must be trivial. So $Y_k\in \langle U\rangle$ is nontrivial for some $k$, that is, $Y_k=_G\Pi u_{t_i}^{s_i}$ with $\sum|s_i|>0$. We still need to show that at least one of them can be expressed in this form with $\sum|s_i|\leq D$.

Pick out all such $k$, and switch all those $u_i$ to $w_i$ to form a new word $z'$, which is also equal to $1$ in $G$. After combining the terms without free reduction, $z'=X'_1Y'_1\dots X'_tY'_t$ with $X'_1$ and $Y'_t$ possibly be empty.

Apply Corollary \ref{R1} again, there exists $k$ such that $X'_k\in \langle W\rangle$ or $Y'_k\in \langle U\rangle$ that is not $1$ as a string of letters. But it can not happen to $Y'_k$ by our construction, hence $X'_k\in \langle W\rangle$ for some $k$.

Now, 
\begin{equation}\label{*}
X'_k=X_jW_jX_{j+1}\dots X_{l-1}W_{l-1}X_{l}=W_0,\end{equation}
where each $W_r$ is a nontrivial combination in $W^{\pm1}$ for $j\leq r\leq l-1$. 
Note that it may start or end with $W_j$, $W_{l-1}$, in which case, we have
$$X_{j+1}W_{j+1}\dots W_{l-2}X_{l-1}=W_j^{-1}W_0W_{l-1}^{-1}.$$ 
So after re-indexing, we may assume that all terms in the middle expression of \eqref{*} are nontrivial. Note that $l>j$.

We prove by contradiction.  Suppose the number of $w_i^{\pm 1}$ in each $W_r$ is greater than $D$ for $j\leq r\leq l-1$, them by Proposition \ref{uncanceled}, the length of reduced form of $W_r$ is also greater than $D$.

First, consider the cancellation between $X_j$ and $W_j$. Suppose the length in $W_j$ that is cancelled by $X_j$ is $L$, we prove that $L< S+T$.

Suppose $W_j=z_1z_2\dots z_h\dots=s_1s_2\dots s_h\dots$, where $s_i$ is the stem of $z_i$ in this combination, and the cancellation ends in the middle or at the end of $s_h$.  Let $z_h=ab$ such that $a$  is canceled by $X_j$ but $b$ is remained. Note that $a$ must be nonempty while $b$ might be an empty word.

Then $X_j=Xa^{-1}z_{h-1}^{-1}\dots z_1^{-1}$ where $X$ is a word such that there is no cancellation between $X$ and $b$. We can rewrite it as $X_j=Xbb^{-1}a^{-1}z_{h-1}^{-1}\dots z_1^{-1}=(Xb)(ab)^{-1}z_{h-1}^{-1}\dots z_1^{-1}$, so $Xb\in X_j\langle W\rangle$. Note that $b$ is the uncanceled part of $s_h$ together with the right twig, or just the right twig, or empty. In whichever case, $b<S+T$.

If $L\geq S+T$, on the one hand,
$$|X_j|=|X|+|a^{-1}z_h^{-1}\dots z_1^{-1}|\geq|X|+S+T.$$

On the other hand,
$$|Xb|=|X|+|b|<|X|+S+T.$$

So $|Xb|<|X_j|$, which contradicts the assumption that $X_j$ is $W$-free.

Similarly, the length of the two segments of $W_i$ canceled by $X_i$ and $X_{i+1}$ are both shorter than $S+T$. The remaining part is at least $D-2(S+T)=3C$.

Next, if $X_i$ is completely canceled, the length that remaining $W_{i-1}$ and $W_i$ cancel each other is no longer than $C$. We again prove it by contradiction.

Suppose the remaining parts are $W_{i-1}'$ and $W_i'$ respectively. The canceled part can be viewed as the common initial part of $W_{i-1}'^{-1}$ and $W_i'$. If the length of this part is longer than $C$, by Lemma \ref{C}, there is a trivial relating pair. Now if we remove this element from both combinations, it separates both words into two parts and the left part is exactly the same.

Now this element, $z$, separates $W_{i-1}$ into two parts, say $A'$ and $A$, so $W_{i-1}=A'zA$. Note that the process is to take one element away from a combination of $W^{\pm1}$, so $A'$ and $A$ are both combinations of elements in $W^{\pm1}$ too. Similarly, it divides $W_{i}$ into $B$ and $B'$, so $W_{i}=Bz^{-1}B'$ with $B$ and $B'$ both all in $\langle W\rangle$. 

Consider the segment $AX_iB$ in $W_{i-1}X_iW_i=A'zAX_iBz^{-1}B'$. From the construction, $A$ and $B$ cancel $X_i$ completely and the remaining parts are inverse of each other, which implies $AX_iB=1$. But this implies $X_i\in \langle W\rangle$, which contradicts the assumption that $X_i$ is a nontrivial  $W$-free word.

Finally, after the cancellation by $X_i$ and $X_{i+1}$ and then by $W_{i-1}$ and $W_{i+1}$, the length of the remaining part of $W_i$ is longer than $D-2(S+T)-2C\geq C$, and denote this part by $\bar{W_i}$. Also denote the remaining part of $X_i$ after free reduction by $\bar{X_i}$, which is well-defined since no $W_i$ is cancelled completely. Then $\bar{X_j}\bar{W_j}\bar{X_{j+1}}\dots\bar{W_{l-1}}\bar{X_l}$ is a reduced word that equals to $W_0$.

Now $\bar{W_j}$ is a segment of $W_j$ that is the same as a segment of $W_0$, so they constitute an identical segment. By Lemma \ref{C}, there is a trivial relating pair. Suppose the corresponding element is $v$ and it divides $W_j$ into $H$ and $H'$, and divides $W_0$ into $R$ and $R'$. Then $X_jH=R$, which contradicts the assumption that $X_j$ is a nontrivial $W$-free word. 

So there exists $j\leq s\leq l-1$ such that $|W_s|\leq D$. Now $W_s$ is a combination of elements in $W^{\pm1}$ and suppose the number of elements contained in it is $N$. by proposition \ref{uncanceled}, $|W_s|\geq N$, which implies $N\leq D$. By our construction, there exists $k$ such that $Y_k=_G\Pi u_{t_i}^{s_i}$ with $0<\Sigma|s_i|\leq D$.
\end{proof}

\subsection{Single Bounded Theorem}

Now we are ready to prove the Single Bounded Theorem, which states that the area of a word can be linearly bounded by the word length in one factor.

\begin{manualtheorem}{\ref{SBT}}
Let $G$ be an amalgam of two free groups $F_1$ and $F_2$ amalgamated along finitely generated subgroups $A_1\in F_1$ and $A_2\in F_2$. If $A_1$ is a malnormal subgroup, then there exists a presentation of $G$ and a constant $B$ such that for any word $z=X_1Y_1\dots X_nY_n=_G1$, $X_i\in F_1$, $Y_i\in F_2$, we have ${\rm{Area}}(z)\leq B\sum^n_{i=1}|X_i|$ in this presentation. 
\end{manualtheorem}

\begin{proof}

Pick an $N$-reduced basis $W=\{w_1,\dots,w_m\}$ for $A_1$. Suppose the corresponding basis in $F_2$ is $U=\{u_1,\dots,u_m\}$ so that $w_i=u_i$ in the amalgam. Then $G$ has a presentation $\langle F_1,F_2;w_i=u_i\rangle$. Since $W$ is N-reduced, we will use previous results to prove it under this presentation.

If $X_1=1$, $Y_1\dots X_nY_n=_G1$ implies $X_2Y_2\dots X_n(Y_nY_1)=_G1$. Similarly, if $Y_n=1$, $X_1Y_1\dots X_n=_G1$ implies $(X_nX_1)Y_1\dots X_{n-1}Y_{n-1}=_G1$. We call this process a conjugation by initial or terminal term. If any other terms are trivial, combine the adjacent terms and $\sum|X_i|$ is not getting bigger in any of these cases. So without loss of generality, assume all terms in the expression of $z$ are nontrivial.

 Applying the following algorithm reduces $z$ to 1.

\emph{Step 1: Each $X_i$ is of the form $W_iX'_iW_i'$, where $X'_i$ is $W$-free (possibly 1) and $W_i$ and $W_i'$ are combinations in $W^{\pm1}$. Switch all  $w_i^{\pm1}$ in $W_i$ and $W_i'$ to $u_i^{\pm1}$. After appropriate combining and free reduction and repeating the process if necessary, and maybe after one conjugation by initial or terminal term, we now have $z'_1=X'_1Y'_1\dots X'_tY'_t$ with all terms nontrivial and all $X'_i$ are $W$-free. The first part of Proof of Lemma~\ref{D} shows that $t\neq 1$. If $t>1$, go to step 2.  If $t=0$, stop.}

\emph{Step 2: According to our construction, all $X'_i$ are $W$-free and nontrivial now. By Lemma \ref{D}, there exists $k$ such that $Y'_k=_GU_k$, where $U_k$ is a combination of no more than $D$ elements in $U^{\pm1}$. Switch all $u_i^{\pm1}$ to $w_i^{\pm1}$ in such $U_k$. After appropriate combining and free reduction and repeating the process if necessary, and maybe after one conjugation by initial or terminal term, we now have $z_1=X^1_1Y^1_1\dots X^{1}_{n_1}Y^{1}_{n_1}$ with all terms nontrivial. If $n_1=0$, stop. If $n_1\neq 0$, go to step 1.}

A cycle in the algorithm is a combination of step 1 and step 2 or step 1 itself if the algorithm stops at step 1.

In step 2, we combined at least two terms. So after finite many cycles, $z$ reduces to $1$ though a sequence $z,z_1,\dots,z_k=1$.

Now we consider the total length of words that are in $F_1$, which was $\sum|X_i|$ initially. In step 1, we rewrite $X_i$ as $W$-free form $W_iX'_iW_i'$. According to the proof of Lemma \ref{D}, the length of cancellation between $X'_i$ and $W_i$ or $W_i'$ is no longer than $S+T$, and if $X'_i$ is totally canceled, the length of cancellation between $W_i$ and $W_i'$ is no longer than $C$. Based on these, we get
\begin{equation}\label{switch}
\begin{split}
|X_i|&=|W_iX'_iW_i'|\\
&\geq |W_i|+|X_i'|+|W_i'|-(4S+4T+2C)\\
&\geq n(W_i)+|X_i'|+n(W_i')-(4S+4T+2C)\\
\implies &|X_i|-|X_i'|\geq n(W_i)+n(W_i')-(4S+4T+2C)
\end{split}
\end{equation}
where $n(W_i)$ denotes the number of elements in $W_i$ when expressed as a reduced combination in $W^{\pm1}$. \eqref{switch} means that in step 1, the length of $i$-th term in $F_1$ is decreased by at least $n(W_i)+n(W_i')-(4S+4T+2C)$.

Now let us estimate how many times we need to rewrite a term in $W$-free form. In first cycle, we did it for every such terms, so at most $n$ times. If this term is not combined with other terms during step 2, we can keep it. We combine two such terms when the term (in $F_2$) between them is switched to a combination in $W^{\pm1}$, and it happens at most $n$ times. So during the whole process, this happens at most $2n$ times. Note that a conjugation by initial or terminal term is also viewed as a combining of terms.

 Suppose in $j$-th cycle, the number of $w_i^{\pm1}$ that are switched to $u_i^{\pm1}$ in step 1 is $r_j$, and the number of $u_i^{\pm1}$ that are switched back to $w_i^{\pm1}$ in step 2 is $t_j$  and the number of terms $U_k$ involved is $s_j$. Note that $t_j\leq Ds_j$. There may be some reductions in each step, so after first cycle, the total length in $F_1$ is bounded by $\sum|X_i|-r_1+R_1(4S+4T+2C)+t_1S$, while $R_1$ is the number of times that we rewrite a term in $F_1$ as $W$-free form. Keep doing it and at last, the length of $z_k=1$ is bounded by $\sum|X_i|-\sum r_j +\sum R_j(4S+4T+2C)+\sum t_jS$ . So we have 
\begin{equation}\label{a}
\begin{split}
0&\leq \sum|X_i|-\sum r_j+\sum R_j(4S+4T+2C)+\sum t_jS\\
&\leq \sum|X_i|-\sum r_j+2n(4S+4T+2C)+DS\sum s_i.
\end{split}
\end{equation}

When we switch a $U_k$ away, one term in $F_2$ is gone. So this can be done at most $n$ times, which implies $\sum s_j\leq n$.

Based on the reduction process, there is a diagram whose area is the total number of times $w_i^{\pm1}$ and $u_i^{\pm1}$ are interchanged, which is $\sum r_j+\sum t_j$. The area of $z$ is defined to be the least area of all such diagrams. So an upper bound for $\sum r_j+\sum t_j$ is also an upper bound for $Area(z)$.

By \eqref{a} and $\sum s_j\leq n$,
\begin{equation}
\begin{split}
\sum r_j&\leq \sum|X_i|+2n(4S+4T+2C)+DSn \\
&\leq (DS+8S+8T+4C+1)\sum|X_i|.
\end{split}
\end{equation}

Combining with $\sum t_j\leq Dn\leq D\sum|X_i|$, the final bound is given by: 
\begin{equation}
\begin{split}
{\rm{Area}}(z)&\leq\sum t_j+\sum r_j\\
&\leq (D+DS+4S+4T+4C+1)\sum |X_i|.
\end{split}
\end{equation} 
\end{proof}

\section{Geometric methods}\label{proof}
In this section, we consider the geometric aspects of the mapping tori. We will explore more about the van Kampen diagrams and some operations on it. By an operation we mean to change a van Kampen diagram to a new one by a certain procedure while the boundary label of the diagram is unchanged.

Without loss of generality, we assume all relators are reduced words, and words in consideration are also reduced, so the boundary label of a van Kampen diagram is always a reduced word.

\subsection{Foldings on diagrams}
Similarly as in geometry, a \emph{path} in a van Kampen diagram is a continuous map $f$ from $[0,N]$ to the one-skeleton of the diagram. Here, we require it to be combinatorial: it is a homeomorphism from $[k,k+1]$ to a $1$-cell for each $k\in \{0,\dots,N-1\}$. Concatenation of paths are defined in the usual way. A path is called a \emph{loop} if $f(0)=f(N)$. A path is called \emph{simple} if it does not contain a backtrack, that is, $f([k-1,k])$ and $f([k,k+1])$ are not inverse edges for any $k$. The \emph{label} of a path is the label of its image. It is \emph{reduced} if the label is a reduced word. In this paper, we do not distinguish a path and its label.

If a path has backtracks, we can remove all backtracks from this path and it results in a simple path. Note that we do not change the diagram here, just change this path. A simple path may not be reduced, We need to change the diagram in order to reduce it.

In last section, Stallings foldings are applied to eliminate unreduced labels in a graph. However for two-dimensional complexes, like planar van-Kampen diagrams, we cannot guarantee all paths are reduced. Yet we still have some operations to fold certain edges based on our needs.

Here we introduce two kinds of foldings on a van Kampen diagram. One is called directed folding and one is called rotated folding.

Suppose there is a path passing through three different vertices $A$, $B$ and $C$ consecutively, where $AB=a$ and $BC=a^{-1}$, we call such a segment $ABC$ an \emph{unreduced segment}. Note that it is not a backtrack. To fold this part, we need to consider whether there are other edges with $B$ as a vertex on two sides of $ABC$.

\begin{figure}
	\includegraphics[width=100mm]{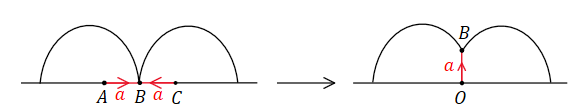}
	\caption{Direct folding}
	\label{fold1}
\end{figure}

\textbf{Direct folding}

If there are no such edges on one side of the path, we can fold it in that side. We call it \emph{direct folding}. Figure \ref{fold1} illustrates this folding. 

After this folding, the unreduced segment turns into a backtrack and we can remove it from the path.

Note that if there are no such edges on the other side either, it results in a vertex of valance 1 on this folded edge, and we call such an edge a \emph{branch}. We can remove this branch by deleting this edge but keep the vertex that also on other edges.

\begin{figure}
	\includegraphics[width=100mm]{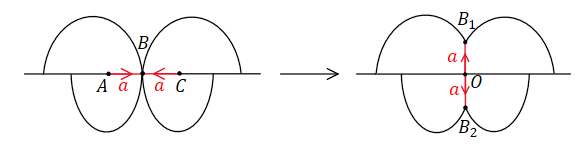}
	\caption{Rotated folding}
	\label{fold2}
\end{figure}

\textbf{Rotated folding}

Suppose there are edges with vertex $B$ on both sides of the segment. Locally, $ABC$ divides the diagram into two parts and it is on the boundary of the closure of each part. So what we do here is do direct folding to each part, and glue them back along the new boundary. The resulted folding is called \emph{rotated folding}. Figure \ref{fold2} illustrates this folding.

After this folding, $A$ and $C$ are merged into one point so we can remove this unreduced segment from the path.

If another segment of this path also passes $B$ from one side of $ABC$ to the other side, we call it a \emph{crossing} at $B$. In this case, we need to insert a segment $a^{-1}a$ to this segment after the folding, so we can not really reduce it on this path. But if a path has no crossings, the following lemma is straightforward.

\begin{lemma}\label{path}
Let $D$ be a van-Kampen diagram and  $\alpha$ is a path on $D$ without crossings. Then there is an operation on $D$, called reducing $\alpha$, such that $\alpha$ is transformed to a reduced path which is the reduced form of $\alpha$, and the area of the diagram is unchanged after this operation..
\end{lemma}

\subsection{The induced map on diagrams}

A group can always be presented by $G=\langle F;R\rangle=F/\llangle R\rrangle$, where $F$ is a free group and $\llangle R\rrangle$ denotes the normal closure of $R$.

If $\Phi$ is an automorphism of $F$ that maps $\llangle R\rrangle$ to $\llangle R\rrangle$, it induces an automorphism $\varphi$ on the quotient group $F/\llangle R\rrangle$. If an automorphism $\varphi$ of $G$ comes from an automorphism $\Phi$ of $F$ in this way, we say $\varphi$ is \emph{induced from the free group automorphism $\Phi$}.

Suppose an automorphism $\varphi$ of $G$ is induced from the free group automorphism $\Phi$. Given a van Kampen diagram $D$ of a word $z=_G1$ in $G$, there is another van Kampen diagram $\Phi(D)$ as follows. 

First consider a cell, which is labeled by a relator $r\in R$. Replace each edge representing a generator by a segment representing its image under $\Phi$. Now there may be some unreduced segments, we can do direct folding to each unreduced segment and delete each branch created, then fill the new loop with minimal number of cells. This process makes sense since $\Phi$ maps $\llangle R\rrangle$ to $\llangle R\rrangle$. The new diagram is a van Kampen diagram with boundary labeled by the reduced form of $\Phi(r)$. If we fix the way we fill each such loop, it defines a map from all cells of $G$ to van Kampen diagrams of $G$.

Now $D$ is composed of cells of $G$ and labeled by generators of $G$. Suppose all the cells are $c_1,\dots,c_n$. First we replace each edge representing a generator by a segment representing its image. Now the boundary of each $c_i$ is replaced by a new loop, and we consider these loops.

If there is an unreduced segment on one loop, we can do direct folding to reduce it, and delete branches if there is any. After finite many times, we can eliminate all unreduced segments on all loops. 

Now all loops are labeled by a reduced word which is the image of the boundary of the corresponding cell.  We can now fill each loop in the fixed way. 
This gives a new van Kampen diagram. If the boundary of this diagram is not reduced, we can do direct foldings to reduce them. After that, we end up with a van Kampen diagram with boundary labeled by reduced form of $\Phi(z)$, and we denote this diagram by $\Phi(D)$. 

Even though $\Phi(D)$ is not unique, it depends on the order that we fold all unreduced segments and the filling of the loops, the area of $\Phi(D)$ is well-defined.

Furthermore, if $R=\{r_1,\dots,r_m\}$ is finite, there are finitely many different cells which we denote by $C_1,\dots,C_m$ respectively. Let $A_i$ denotes the area of $\Phi(C_i)$. Let $A=\max A_i$, then for any van Kampen diagram $D$, ${\rm{Area}}(\Phi(D))\leq A{\rm{Area}}(D)$.

A special case is when $A=1$, then the image of one cell is exactly another cell. In this case, the image of a van Kampen diagram has the same area and we say that $\Phi$ is \emph{area preserving}. 

\subsection{Room moving}
For a van Kampen diagram of the mapping torus $M_{\varphi}=\{G,t;t^{-1}at=\varphi (a), a\in G\}$, the cells include those from relations in $G$ and those from the relations $t^{-1}at=\varphi (a)$. We call the cells of the first type \emph{primitive cells}, and the cells of the second type \emph{$t$-cells}.

\begin{figure}
\includegraphics[width=90mm]{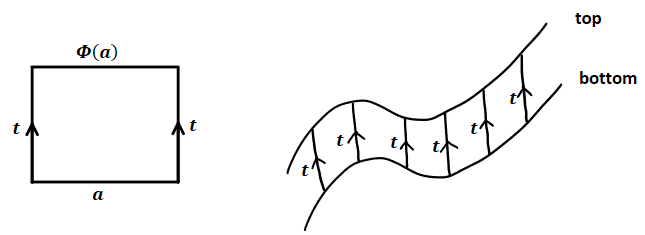}
\caption{}
\label{t-corridor}
\end{figure}

In the diagrams of mapping tori, $t$-corridors have been widely studied, see \cite{MD}. Simply speaking, a $t$-corridor is a series of $t$-cells connected along $t$-sides. The Figure \ref{t-corridor} is an example of a $t$-cell and a $t$-corridor.

In this paper, $t$-corridors are assumed to be maximal. That is, we cannot extend it to a longer one from either side. When it is clear from the context, we also call a $t$-corridor simply a corridor.

Depending on the structure of the $t$-cell, a $t$-corridor either starts and ends on the boundary, or forms a ring which we call a \emph{$t$-ring}. Due to the following lemma, we only need to consider corridors of the first type.

\begin{lemma}
In a van Kampen diagram, there is an operation that can remove each $t$-ring. If $\Phi$ is area preserving, the number of primitive cells is unchanged after this operation.
\end{lemma}
\begin{proof}
\begin{figure}
	\includegraphics[width=90mm]{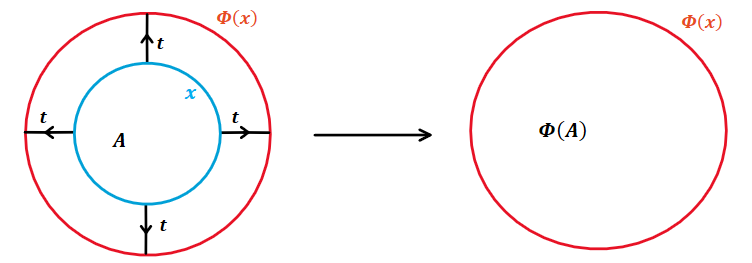}
	\caption{Remove a t-ring}
	\label{t-ring}
\end{figure}

The tops and bottoms of all $t$-rings do not have crossings with each other, so first we reduce them simultaneously.

Suppose the inner loop in a $t$-ring is labeled by $x$, and the enclosed subdiagram is $A$. We can remove everything inside the outer loop, including $A$ and the $t$-ring, and refill it with $\Phi(A)$. Figure \ref{t-ring} illustrates this operation.
 
If the $t$-sides are facing inward in a $t$-ring, replace $\Phi$ by $\Phi^{-1}$ and the same result follows.

If $\Phi$ is area preserving, So is $\Phi^{-1}$, and the area of $\Phi(A)$ or $\Phi^{-1}(A)$ is the same as area of $A$, so the number of primitive cells is unchanged after each operation.
\end{proof}

For the rest of this paper, all corridors are assumed to start and end on the boundary.
In this case, a corridor separates the diagram into two parts, one is on the bottom side and one is on the top side. So we can talk about relative position of other cells or subdiagrams and this corridor.

\begin{definition}
In a van Kampen diagram, a room is the closure of one component after deleting all $t$-corridors. A $t$-corridor \emph{bounds a room} if they have non-empty intersection.
\end{definition}

Note that if a $t$-corridor bounds a room, their intersection is a connected path or a point.

Here we introduce another kind of operation on a diagram called \emph{room moving}. That is, we can essentially move a room across a corridor.

\begin{lemma}\label{7}
If a $t$-corridor bounds a room in a van Kampen diagram, we can do an operation that moves all primitive cells in this room to the other side of the corridor. If $\Phi$ is area preserving, the number of primitive cells in the diagram is unchanged after the operation.
\end{lemma}

\begin{proof}

\begin{figure}
	\includegraphics[width=120mm]{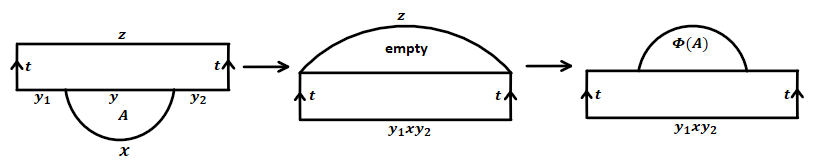}
	\caption{}
	\label{move1}
\end{figure}

\begin{figure}
	\includegraphics[width=120mm]{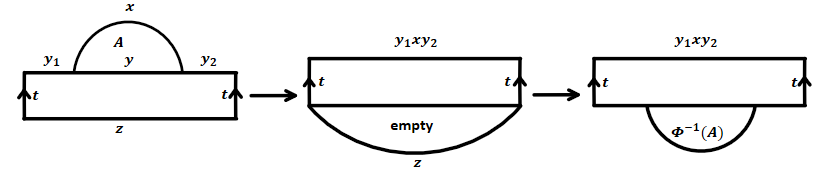}
	\caption{}
	\label{move2}
\end{figure}

Suppose the room $A$ intersects the corridor with a path $y$, and $A=xy^{-1}$ with $x$ is a segment of boundary of $A$ that not on the corridor. We can choose $x$ and $y$ so that the loop $xy^{-1}$ has no crossings. Note that $x$ or $y$ may be a point. 

If $x$ is not a simple path, we can get a simple path by deleting all backtracks. Note that all cells in $A$ are between this simple path and $y$, so we can replace $x$ by this simple path. Without loss of generality, we assume that $x$ is a simple path.

First we consider the case that $y$ is on the bottom of the corridor. Let the bottom be $y_1yy_2$ and the top be $z$, then $z=\Phi(y_1yy_2)=_G\Phi(y_1xy_2)$. Now we consider the subdiagram containing this corridor and the room $A$. First by Lemma \ref{path}, we can reduce its boundary. So we assume its boundary is already reduced, that is, $y_1xy_2$ and $z$ are both reduced words.

What we will do here is remove other parts of the diagram, reconstruct this part, then attach other parts back in the same way.

First, construct a $t$-corridor with bottom $y_1xy_2$ and top the reduced form of $\Phi(y_1xy_2)$. Second, attach a path $z$ to two ends of the top to form a loop $\Phi(y_1xy_2)z^{-1}$ with interior disjoint from the corridor. Note that $\Phi(y_1xy_2)z^{-1}=\Phi(xy^{-1})$, which means this boundary is the image of the boundary of $A$. Third, fold two ends of $z$ with top of the corridor, if any, so the new simple loop is reduced. Forth, fill this loop with $\Phi(A)$. Now the boundary is still $y_1xy_2tz^{-1}t^{-1}$, so we can attach other parts of the diagram back.

The process is similar if $y$ is on the top of the corridor. Let the top be $y_1yy_2$ and the bottom be $z$, then $y_1yy_2=\Phi(z)$. we need to construct a $t$-corridor with bottom the reduced form of $\Phi^{-1}(y_1xy_2)$ and top $y_1xy_2$, which is assumed to be reduced now. Then attach a path $z$ to the bottom and fold two ends, fill the loop with $\Phi^{-1}(A)$.

These two cases are shown in Figure \ref{move1} and Figure \ref{move2}, where the first step is constructing the $t$-corridor and attaching the path $z$, the second step is folding it and filling the new room.
	
In both cases, we replace $A$ by its image under $\Phi$ or $\Phi^{-1}$, which is now on the different side of the corridor. So we can view this process as moving $A$ to the other side of the corridor. Except $A$ and the corridor, any other parts of the diagram are unaffected. So if $\Phi$ is area preserving, the total number of primitive cells is unchanged.
\end{proof}

\subsection{Bound for primitive cells}
To prove Theorem \ref{1}, we need to evaluate the number of cells in a van Kampen diagram, including primitive cells and $t$-cells. In this section, we give a bound for the primitive cells.

First we fix a presentation for $G$: pick an N-reduced basis $W=\{w_1,\dots,w_m\}$ for $A_1$, and suppose the corresponding basis in $F_2$ is $U=\{u_1,\dots,u_m\}$ so that the presentation is $\langle F_1, F_2; w_1=u_1,\dots,w_m=u_m \rangle$.

\begin{lemma}\label{primitive bound}
Let $z=X_1Y_1t^{\lambda_1}X_2Y_2t^{\lambda_2}\dots X_nY_nt^{\lambda_n}$ be a word that equals 1 in $M_\varphi$, $X_i\in F_1, Y_i\in F_2$. There exists a van Kampen diagram with number of primitive cells bounded by $B\sum|X_i|$, where $B$ is the constant from the Single Bounded Theorem \ref{SBT}.
\end{lemma}

\begin{proof}
In the mapping torus, we have $\varphi(x)=t^{-1}xt$, which implies $t^{-1}x=\varphi(x)t^{-1}$, and $tx=\varphi^{-1}(x)t$. We can use these relations to push $t$ to the end of the expression. Since $z=_{M_{\varphi}}1$, $\Sigma \lambda_i=0$. We have a series of transformations:
\begin{equation}
\begin{split}
z&=X_1Y_1t^{\lambda_1}X_2Y_2t^{\lambda_2}\dots X_nY_nt^{\lambda_n}\\
&\rightarrow X_1Y_1\varphi^{-\lambda_1}(X_2)\varphi^{-\lambda_1}(Y_2)t^{\lambda_1}t^{\lambda_2}\dots X_nY_nt^{\lambda_n}\\
&\rightarrow \dots\dots\\
&\rightarrow X_1Y_1\varphi^{-\lambda_1}(X_2)\varphi^{-\lambda_1}(Y_2)\varphi^{-\lambda_1-\lambda_2}(X_3)\varphi^{-\lambda_1-\lambda_2}(Y_3) \dots \varphi^{-\lambda_1-\dots-\lambda_{n-1}}(X_n)\\
&\ \ \ \ \ \varphi^{-\lambda_1-\dots-\lambda_{n-1}}(Y_n)t^{\lambda_1+\dots+\lambda_n}\\
&\rightarrow X_1Y_1X_2\varphi^{-\lambda_1}(Y_2)X_3\varphi^{-\lambda_1-\lambda_2}(Y_3)\dots X_n\varphi^{-\lambda_1-\dots-\lambda_{n-1}}(Y_n)=z'
\end{split}
\end{equation}

Finally, we can transform $z'$ to 1 with no more than $B\sum|X_i|$ relations by Theorem~\ref{SBT} .

Now, the whole process gives a diagram denoted by $D$. From $z$ to $z'$, we only used relation $t^{-1}xt=\varphi(x)$, thus the corresponding cells are all $t$-cells. The last step is done in $G$, thus the corresponding cells are all primitive cells. So the total number of primitive cells in $D$ is bounded by $B\sum|X_i|$.
\end{proof}

\subsection{Proof of Theorem \ref{1}}
For reader's convenience, we recall the statement of the theorem.

Let $G$ be an amalgam of two finitely generated free groups $F_1$ and $F_2$ amalgamated along finitely generated subgroups $A_i\in F_i$, where $A_1$ is a malnormal subgroup. Let $\psi$ be an automorphism of $F_2$ that fixes $A_2$ pointwise. Then there exists an induced automorphism $\varphi$ of $G$ such that it restricts to identity on $F_1$, and $\psi$ on $F_2$.

\begin{manualtheorem}{\ref{1}}
The Dehn function of the mapping torus $$M_{\varphi}=\{G,t;t^{-1}at=\varphi (a), a\in G\}$$ is quadratic.
\end{manualtheorem}

\begin{proof}
To determine the Dehn function of a group, the most common way is to find both a lower bound and an upper bound. For the lower bound, it is enough to show that the mapping torus is not hyperbolic, which implies its Dehn function is not linear. By the result of Gromov \cite{M}, it is at least quadratic.

Pick any generator $a$ from $F_1$, since $\varphi$ is an identity on $F_1$, we have $t^{-1}at=a$, which shows $a$ and $t$ commute. Furthermore, the inclusion of $F_1$ and $\langle t\rangle$ in $M_{\varphi}$ are both embeddings, so $a$ and $t$ generate a subgroup isomorphic to $\mathbb{Z}^2$, which implies $M_{\varphi}$ is not hyperbolic.

Next, we give a quadratic upper bound for the mapping tours. This will complete our proof.

Based on Lemma \ref{primitive bound}, there is a van-Kampen diagram $D$ with the number of primitive cells bounded by $B\sum|X_i|\leq B|z|$.

Note that $\varphi$ is induced from $\Phi$: $F_1*F_2\rightarrow F_1*F_2$, which is an identity on $F_1$, and $\psi$ on $F_2$. It is also area-preserving since it maps each cell to itself.

\begin{figure}
\includegraphics[width=110mm]{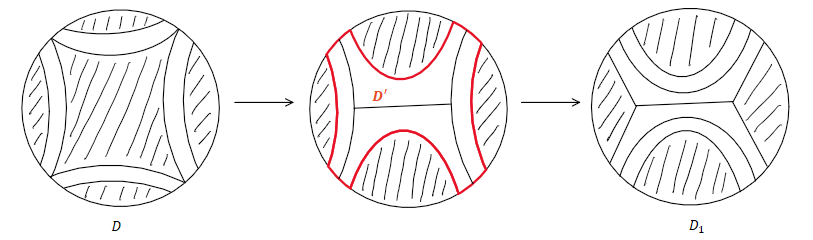}
\caption{}
\label{Transformation}
\end{figure}

First, by the operation from Lemma \ref{7}, we can move all primitive cells towards the boundary so that the set of all $t$-corridors forms a connected, simply-connected subdiagram, which is denoted by $D'$. The boundary of $D'$ are either boundary of some primitive cells or boundary of the original diagram. By Lemma \ref{7}, the number of primitive cells does not change after the operations. Suppose the maximal length of the boundary of each primitive cell is $m$, then the sum of length of boundaries of all primitive cells is bounded by $Bm|z|$. Therefore the length of boundary of $D'$ is bounded by $Bm|z|+|z|=(Bm+1)|z|$.

The next step is refilling for $D'$. Note that the map is induced from $\Phi$, which is an automorphism of the free group $F_1*F_2$, so $D'$ can be viewed as a van-Kampen diagram of the mapping torus of this free group automorphism. According to Theorem \ref{3}, there is a refilling of $D'$ with $t$-cells whose number is bounded by $M(Bm+1)^2|z|^2$ for some constant $M$ depending on $\Phi$. This gives a new diagram, denoted by $D_1$. Figure \ref{Transformation} is an example of this series of operations, in which we shaded all primitive cells.

Since ${\rm{Area}}(z)$ is defined as the least area of all such diagrams, the bound for area of $D_1$ also gives an upper bound for ${\rm{Area}}(z)$:
$${\rm{Area}}(z)\leq {\rm{Area}}(D_1)\leq M(Bm+1)^2|z|^2+B|z|.$$
\end{proof}

\end{document}